\documentclass[12pt,reqno]{amsart}
\usepackage{enumerate, latexsym, amsmath, amsfonts, amssymb, amsthm, graphicx, color,bm}
\usepackage{hyperref}
\hypersetup{hypertex=true,
	colorlinks=true,
	linkcolor=blue,
	anchorcolor=blue,
	citecolor=blue}
\usepackage{amsmath}
\def\pmod #1{\ ({\rm{mod}}\ #1)}
\def\Z{\Bbb Z}

\def\Q{\Bbb Q}

\def\bg{\bigg}
\def\({\bg(}
\def\){\bg)}

\def\adj{{\rm adj}}

\def\u{{\bm u}}
\def\v{{\bm v}}

\def\Ack{\medskip\noindent {\bf Acknowledgments}}
\theoremstyle{plain}
\newtheorem{theorem}{Theorem}

\newtheorem{lemma}{Lemma}

\newtheorem{conjecture}{Conjecture}
\theoremstyle{definition}

\theoremstyle{remark}
\newtheorem{remark}{Remark}

\makeatletter
\@namedef{subjclassname@2010}{%
	\textup{2010} Mathematics Subject Classification}
\makeatother
\vspace{4mm}

\allowdisplaybreaks[1]
\begin{document}

	\title{On the cyclotomic field $\mathbb{Q}(e^{2\pi {\bf i}/p})$ and Zhi-Wei Sun's conjecture on $\det M_p$}
	
	\author[L.-Y. Wang, H.-L. Wu]{Li-Yuan Wang and Hai-Liang Wu*}
	
	\address {(Li-Yuan Wang) School of Physical and Mathematical Sciences, Nanjing Tech University, Nanjing 211816, People's Republic of China}
	\email{\tt wly@smail.nju.edu.cn}
	
	\address {(Hai-Liang Wu) School of Science, Nanjing University of Posts and Telecommunications, Nanjing 210023, People's Republic of China}
	\email{\tt whl.math@smail.nju.edu.cn}

	\begin{abstract}
		In 2019, Zhi-Wei Sun posed an interesting conjecture on certain determinants with Legendre symbol entries. In this paper, by using the arithmetic properties of $p$-th cyclotomic field and the finite field $\mathbb{F}_p$, we confirm this conjecture. 
	\end{abstract}
	
	\thanks{
		\newline \indent
		This research was supported by the National Natural Science Foundation of China (Grant Nos. 12201291 and 12101321). 
		\newline\thanks{*Corresponding author.}}
	\maketitle
	\section{Introduction}	
	\setcounter{lemma}{0}
	\setcounter{theorem}{0}
	\setcounter{corollary}{0}
	\setcounter{remark}{0}
	\setcounter{equation}{0}
	\setcounter{conjecture}{0}
	\subsection{Notations} 
	We first introduce the main notations we shall use. Let $p$ be an odd prime and $\mathbb{F}_p$ be the finite field with $p$ elements. Let $(\frac{\cdot}{p})$ denote the Legendre symbol, i.e., 
	$$
	\left(\frac{x}{p}\right)=\begin{cases}
		0  & \mbox{if}\ x=0,\\
		1  & \mbox{if}\ x\ \text{is a nonzero square over $\mathbb{F}_p$},\\
		-1 & \mbox{otherwise.}
	\end{cases}$$
    Throughout this paper, we let  $\zeta =e^{2\pi {\bf i}/p}$, where ${\bf i}\in\mathbb{C}$ is the primitive $4$-th root of unity with argument $\pi/2$. Let $\mathfrak{p}=(1-\zeta)\mathbb{Z}[\zeta]$ be the unique prime ideal of $\mathbb{Z}[\zeta]$ containing $p$. 

	Also, For any matrix $A=[a_{ij}]_{0\le i,j\le m }$, we use  $\det A$ to denote its determinant and $\adj (A)=[A_{ji}]_{0\le i,j\le m}$ to denote its adjugate matrix.
	
	\subsection{Background and Motivation}
	The research on determinants with Legendre symbol entries can be traced back to the work of D. H. Lehmer \cite{Lehmer} and L. Carlitz \cite{carlitz}. For example, Carlitz \cite[Theorem 4 (4.9)]{carlitz} first studied the arithmetic properties of the matrix 
	$$\left[\left(\frac{j-i}{p}\right)\right]_{1\le i,j\le p-1}$$
	and proved that its characteristic polynomial is 
	$$\left(t^2-(-1)^{\frac{p-1}{2}}p\right)^{\frac{p-3}{2}}\left(t^2-(-1)^{\frac{p-1}{2}}\right).$$
	
	Along this line, in 2004 R. Chapman \cite{chapman} investigated many interesting variants of Carlitz's matrices and the most well-known variant is known as the matrix 
	\begin{equation}\label{Eq. definition of Cp}
		C_p:=\left[\left(\frac{j-i}{p}\right)\right]_{0\le i,j\le \frac{p-1}{2}}.
	\end{equation}
	Due to the complexity of the evaluation of $\det C_p$, Chapman \cite{problems,evil} posed the following conjecture which is known as the Chapman's ``evil determinant" conjecture. 
	\begin{conjecture} Let $p$ be an odd prime. Then 
		$$\det C_p=\begin{cases}
			-a_p & \text{{\rm if}}\ p\equiv 1\pmod4,\\
			1    & \mbox{{\rm if}}\ p\equiv 3\pmod4. \\
			
		\end{cases}
		$$
		The number $a_p$ is defined by the following equality
		$$\varepsilon_p^{(2-(\frac{2}{p}))h(p)}=a_p+b_p\sqrt{p}\ (a_p,b_p\in\mathbb{Q}),$$
		where $\varepsilon_p>1$ and $h(p)$ denote the fundamental unit and class number
		of $\mathbb{Q}(\sqrt{p})$ respectively. 
	\end{conjecture}
	By using sophisticated matrix decomposition techniques, M. Vsemirnov \cite{V1,V2} confirmed this conjecture completely. 
	
	Z.-W. Sun \cite{S19} further studied some variants of the above results. Let $M_p$ be the matrix obtained from $[(\frac{i-j}{p})]_{0 \le i,j \le (p-1)/2}$ via replacing all the entries in the first row by $1$, i.e., 
	$M_p=[e_{ij}]_{0\le i,j\le (p-1)/2}$ with
	$$e_{ij}=\begin{cases}1&\mbox{if}\ i=0,\\ (\frac{i-j}{p})&\mbox{otherwise}.\end{cases}$$
	Sun \cite[Conjecture 4.6(i)]{S19} posed a conjecture involving the explicit value of $\det M_p$. 
	
	\subsection{Main Results}
	
In this paper, we confirm Sun's conjecture on $\det M_p$ completely. Namely, we have the following result.
	
	\begin{theorem}\label{th} Let $p$ be an odd prime. Then 
		$$\det M_p=
		\begin{cases}
			(-1)^{\frac{p-1}{4}}&\text{{\rm if}}\ p \equiv 1 \pmod 4,\\
			(-1)^{(h(-p)-1)/2}&\mbox{{\rm if}} \ p\equiv 3 \pmod 4,
		\end{cases}
		$$	
		where $h(-p)$ is the class number of $\Q(\sqrt{-p})$.
	\end{theorem}
	
	\begin{remark}
		Our idea to  prove  Theorem \ref{th} is as follows. We first show that  $\det M_p\in \mathbb{Z}[\zeta]^\times$
		(where $\mathbb{Z}[\zeta ]^\times$ denotes the unit group of  $\mathbb{Z}[\zeta]$), i.e., there is an algebraic integer $y\in\mathbb{Z}[\zeta]$ such that 
		$$y\cdot \det M_p=1.$$
		As $\det M_p\in\mathbb{Z}$, we see that the algebraic integer $y$ is indeed a rational number and hence $y\in\mathbb{Z}$. By this, the equality $y\cdot \det M_p=1$ implies that $\det M_p\in\{\pm1\}$. 
		
		We next prove that Theorem \ref{th} is true in the sense of modulo $p$, i.e., 
		\begin{align}\label{mod}
			\det M_p\equiv
			\begin{cases}
				(-1)^{\frac{p-1}{4}}\qquad\;\; \pmod p&\mbox{if}\ p \equiv 1 \pmod 4,\\
				(-1)^{(h(-p)-1)/2}\pmod p &\mbox{if} \ p\equiv 3 \pmod 4.
			\end{cases}
		\end{align}
		Combining the above results, Theorem \ref{th} follows directly. 
	\end{remark}
	
	\subsection{Outline of the paper} We will divide our proof of Theorem \ref{th} into two parts. In Section 2, we shall show that $\det M_p$ is a unit of $\Z[\zeta]$. The congruence (\ref{mod}) will be proved in Section 3. 
	\maketitle

	\section{First part of the proof :  $\det M_p\in \mathbb{Z}[\zeta ]^\times$}
	\setcounter{lemma}{0}
	\setcounter{theorem}{0}
	\setcounter{corollary}{0}
	\setcounter{remark}{0}
	\setcounter{equation}{0}
	\setcounter{conjecture}{0}
	
	Recall that $\zeta =e^{2\pi {\bf i}/p}$. In this section we show that $\det M_p$ is a unit in $\mathbb{Z}[\zeta ]$. We begin with the following known result (see \cite[(3.1)]{V2}).
	\begin{lemma}\label{Lem. to prove det N}
		Let $u_i,v_j(1\le i,j\le m)$ be complex numbers such that $u_iv_j\ne -1$ for any $1\le i,j\le m$.
		Then 
		\begin{equation}
			\det \bigg[\frac{1}{1+u_iv_j}\bigg]_{1\le i,j\le m}=\prod_{1\le i<j\le m}^{}\(   (u_i-u_j)(v_j-v_i) \)\cdot \prod_{1\le i,j\le m}^{} (1+u_iv_j)^{-1}.  \nonumber
		\end{equation}
	\end{lemma}
	
	Our method used in this section is inspired by Vsemirnov’s decomposition of $C_p$ defined by (\ref{Eq. definition of Cp}). For convenience, we also follow Vsemirnov's notations \cite{V1,V2} and let $n=(p-1)/2$ throughout the proof below. 
	
	Now we are in a position to prove that $\det M_p\in \mathbb{Z}[\zeta ]^\times$. 
	
	{\noindent\bf First part of the proof of Theorem \ref{th}.} Define square matrices $U, V, D$ of order $n+1$ with
	$(i,j)$ entry as follows:
	\begin{align}
		u_{ij}=&\frac{(\frac{i}{p})\zeta^{-j-2i}+(\frac{-j}{p})\zeta^{-2j-i}}  { \zeta^{-i-j}+(\frac{i}{p})(\frac{-j}{p}) } , \quad     0\le i,j\le n,\\
		v_{ij}=&\zeta^{2ij} ,\quad   0\le i,j\le n ,\\
		d_{ii}=&\prod_{0\le j\le n \atop j\ne i}^{}\frac{1}{\zeta^{2i}-\zeta^{2j}}, \quad    0\le i\le n,\\
		d_{ij}=&0,  \quad    i\ne j.
	\end{align}
	By Vsemirnov's result 
	$$C_p=\(\frac{-1}{p}\)\tau_p(2) \zeta^{(p^2-1)/4}\cdot VDUDV,$$
	where 
	$$\tau_p(a)=\sum^{p-1}_{k=1}\left(\frac{k}{p}\right)\zeta^{ak}$$
	is the Gauss sum. Also, $\tau_p(1)$ will be abbreviated as $\tau$. By the Galois theory 
	$$\tau_p(2)=\(\frac{2}{p}\)\tau=\(\frac{2}{p}\)\sqrt{(-1)^{{(p-1)/2}}p}.$$
	
	Recall that $M_p$ is obtained from $[(\frac{i-j}{p})]_{0 \le i,j \le (p-1)/2}$ via replacing all the entries in the first row by $1$. 
	Define $M_p^-$ to be the matirx obtained from $[(\frac{j-i}{p})]_{0 \le i,j \le (p-1)/2}$ via replacing all the entries in the first row by $1$.  It is easy to check that 
	$$\det M_p=
	\begin{cases}
		\det M_p^-&\mbox{if}\ p \equiv 1 \pmod 4,\\
		-\det M_p^-&\mbox{if} \ p\equiv 3 \pmod 4.
	\end{cases}
	$$	
	Hence  $\det M_p\in \mathbb{Z}[\zeta ]^\times$ is equivalent to  $\det M_p^-\in \mathbb{Z}[\zeta ]^\times$.

 Let  $\adj (C_p)=[C_{ji}]_{0\le i,j\le n} $ be the adjugate matrix of $C_p$. Then, by  expanding $\det M_p^-$ along the first row, we see that  
	\begin{align}
		\det M_p^-=\sum_{k=0}^{n}C_{0k}&=(1,1,\ldots,1)\cdot [C_{ji}]_{0\le i,j\le n}\cdot (1,0,0,\ldots,0)^T   \nonumber  \\
		&=(1,1,\ldots,1)\cdot\adj (C_p)\cdot (1,0,0,\ldots,0)^T. \nonumber
	\end{align}
	Let $\lambda =(\frac{-2}{p})\tau \zeta^{(p^2-1)/4}$. Then 
	
	\begin{align}
		\adj (C_p)= \lambda^n\cdot 	\adj (V) \cdot	\adj (D) \cdot	\adj  (U) \cdot	\adj (VD).\nonumber
	\end{align}
	Hence
	\begin{align}\label{vdudv}
		\det M_p^{-}=\lambda^n\cdot (1,1,\ldots,1)\cdot 	\adj (V) \cdot	\adj (D) \cdot	\adj  (U) \cdot	\adj (VD)\cdot  (1,0,\ldots,0)^T.
	\end{align}
		For any matrix $A=[a_{ij}]_{0\le i,j\le m }$, let $A^{(k)}$ denote the matrix obtained from $A$ via replacing all entries in the $(k+1)$-th row by $1$. One can verify that 
	\begin{align}
		(1,1,\ldots,1)\cdot 	\adj (V) &= (1,1,\ldots,1)\cdot [V_{ji}]_{0\le i,j\le n}       \nonumber\\
		&=\(\sum_{k=0}^{n}V_{0k},\sum_{k=0}^{n}V_{1k},\ldots ,\sum_{k=0}^{n}V_{nk}\  )    \nonumber  \\
		&=\(  \det V^{(0)}  ,\det V^{(1)}  ,\ldots ,\det V^{(n)}\). \nonumber
	\end{align}
	Since each element of the first row of $V$ is $1$, we see that  $\det V^{(0)}=\det V  $ and $ \det V^{(k)}  =0$ for each $1\le k\le n$. Thus
	$
	(1,1,\ldots,1)\cdot 	\adj (V)  = \det V \cdot ( 1  ,  0,0,\ldots,0  ).
	$
	Noting that $D$ is a diagonal matrix with $(i+1)$-th element $d_{ii}$, we have 
	\begin{displaymath}
		\adj (D)=\det D \cdot D^{-1}
		=\det D\cdot 	\left( \begin{array}{cccccccc}
			\frac{1}{d_{00}} & & \  &  \\
			&	\frac{1}{d_{11}} &    &  \\
			& & \ddots & \\
			& &  &	\frac{1}{d_{nn}} \\
		\end{array} \right).
	\end{displaymath}
	Hence 
	$$(1,1,\ldots,1)\cdot 	\adj (V) \cdot	\adj (D)=\det V \cdot \det D \cdot 	\frac{1}{d_{00}} \cdot ( 1  ,  0,0,\ldots,0).$$
	This implies 
	\begin{align}\label{vdu}
		(1,1,\ldots,1)\cdot 	\adj (V) \cdot	\adj (D) \cdot	\adj (U)=\frac{\det V \cdot \det D}{d_{00}}\cdot (U_{00}, U_{10},\ldots,U_{n0}).
	\end{align}
	
	We next consider 
	$\adj (VD)\cdot  (1,0,0,\ldots,0)^T$.
	Let $F=VD$. Then $\adj (F)=[F_{ji}]_{0\le i,j\le n}$ 
	and 
	\begin{align}\label{dv}
		\adj (VD)\cdot  (1,0,0,\ldots,0)^T&=[F_{ji}]_{0\le i,j\le n}\cdot  (1,0,0,\ldots,0)^T      \\
		&=(F_{00}, F_{01},\ldots,F_{0n})^T.\nonumber
	\end{align}
	By (\ref{vdudv})--(\ref{dv}) we obtain 
	\begin{align}
		\det M_p^-=\lambda^n\cdot\frac{\det V \cdot \det D}{d_{00}}\cdot \sum_{k=0}^{n}F_{0k}\cdot U_{k0}.
	\end{align}
	We now evaluate $F_{0k}$. Observing that 
	$$F=[\zeta^{2ij}]_{0\le i,j\le \frac{p-1}{2}}\cdot [d_{ij}]_{0\le i,j\le \frac{p-1}{2}}=[\zeta^{2ij}\cdot d_{jj}]_{0\le i,j\le \frac{p-1}{2}},$$
	one can verify that 
	\begin{align}
		F_{0k}&=(-1)^{0+k}\cdot\det [\zeta^{2ij}\cdot d_{jj}]_{1\le i\le n\atop  0\le j\le n, j\ne k}     \nonumber \\
		&=(-1)^{k}\cdot  \prod_{j=0\atop j\ne k}^{n}d_{jj}\cdot \det  [\zeta^{2ij}]_{1\le i\le n\atop  0\le j\le n, j\ne k}    \nonumber\\
		&=(-1)^{k}\cdot  \prod_{j=0\atop j\ne k}^{n}d_{jj}\cdot  \prod_{j=0\atop j\ne k}^{n} \zeta^{2j}  \cdot   \det  [\zeta^{2(i-1)j}]_{1\le i\le n\atop  0\le j\le n, j\ne k}   \nonumber \\
		&=(-1)^{k}\cdot \frac{ \det D}{d_{kk}}\cdot \zeta^{n(n+1)}  \cdot  \zeta^{-2k}\cdot V(1,\zeta^{2},\ldots,\zeta^{2(k-1)},\zeta^{2(k+1)},\ldots ,\zeta^{2n})      \nonumber \\
		&=(-1)^{k}\cdot \frac{ \det D}{d_{kk}}\cdot \zeta^{n(n+1)}  \cdot  \zeta^{-2k} \cdot \frac{\prod_{0\le i<j \le n}(\zeta^{2j}-\zeta^{2i})}{\prod_{0\le j\le n\atop j\ne k}(\zeta^{2k}-\zeta^{2j})}\cdot (-1)^{n-k}, \nonumber\\\nonumber
	\end{align}
	where $V(a_1,a_2,\ldots,a_m)$ denotes the determinant of the Vandermonde matrix $[a_i^{j-1}]_{1\le i,j\le m}$.
	Noting that 
	$$d_{kk} \cdot \prod_{0\le j\le n\atop j\ne k}(\zeta^{2k}-\zeta^{2j})=1$$
	and 
	$$\det V=\prod_{0\le i<j \le n}(\zeta^{2j}-\zeta^{2i}),$$
	we obtain 
	\begin{align}
		F_{0k}=(-1)^n\cdot  \zeta^{n(n+1)}\cdot \det D\cdot \det V\cdot \zeta^{-2k}.\nonumber
	\end{align}
	Thus 
	$$\det M_p^-=\frac{(-\lambda)^n }{d_{00}} \cdot  \zeta^{n(n+1)}\cdot (\det D)^2 \cdot (\det V)^2\cdot\sum_{k=0}^{n} \zeta^{-2k}\cdot U_{k0}. $$
	Let $H$ be the matrix obtained from $U$ via replacing the first column by $(1,\zeta^{-2},\ldots,\zeta^{-2n})^T$.
	Then 
	\begin{equation}
		\det H=\sum_{k=0}^{n} \zeta^{-2k}\cdot U_{k0}, \nonumber
	\end{equation}
	and 
	\begin{align}
		\det M_p^-&=\frac{(-\lambda)^n }{d_{00}} \cdot  \zeta^{n(n+1)}\cdot (\det D)^2 \cdot (\det V)^2\cdot \det H  \nonumber   \\
		&=\(\frac{2}{p}\)^n\cdot \zeta^{n(n+1)^2}\cdot \tau^n\cdot \prod_{k=1}^{n}(1-\zeta^{2k})\cdot (\det D)^2 \cdot (\det V)^2\cdot \det H  .\nonumber
	\end{align}
	
	For any $1\le k\le p-1$, clearly $\frac{1-\zeta^{k}}{1-\zeta}$
	and 
	$1+\zeta^k$
	are units in $\mathbb{Z}[\zeta]$. Recall that  $\mathfrak{p}=(1-\zeta)\mathbb{Z}[\zeta]$ is the ideal generated by $1-\zeta$.  Then  $$\tau^{2}\mathbb{Z}[\zeta]=p\mathbb{Z}[\zeta]=\mathfrak{p}^{p-1}=(1-\zeta)^{2n}\mathbb{Z}[\zeta].$$ 
	By the unique factorization of nonzero ideals of $\mathbb{Z}[\zeta]$ we see that 
	\begin{align}\label{mu}
		\prod_{k=1}^{n}(1-\zeta^{2k})=\mu \cdot \tau
	\end{align}
	for some $\mu \in \mathbb{Z}[\zeta]^\times$. Also, by Vsemirnov's result \cite{V1,V2} we have 
	$$\tau^{n+1}\cdot (\det D)^2 \cdot (\det V)^2\in\mathbb{Z}[\zeta]^{\times}. $$
	
	In view of the above, to complete the proof of $\det M_p\in\mathbb{Z}[\zeta]^{\times}$, we only need to show that $\det H$ is also a  unit in $\mathbb{Z}[\zeta]$. Define
	\begin{displaymath}
		G =	\left( \begin{array}{cccccccc}
			1& & \  &  \\
			&	 (\frac{1}{p})\zeta^{} &    &  \\
			& &	 (\frac{2}{p})\zeta^{2} &    &  \\
			&	& & \ddots & \\
			&	& &  &	 (\frac{n}{p})\zeta^{n}\\
		\end{array} \right).
	\end{displaymath}
	Then one can verify that 
	\begin{displaymath}
		GHG =	\left( \begin{array}{cccccccc}
			1                                           &  (\frac{-1}{p})                     & (\frac{-1}{p})    &\ldots   &  (\frac{-1}{p}) \\
			(\frac{1}{p})\zeta^{-1}	                 &	 \widetilde{u}_{11}                 &   \widetilde{u}_{12}     &\ldots  &   \widetilde{u}_{1n}\\
			(\frac{2}{p})\zeta^{-2}    &  \widetilde{u}_{21}                     & \widetilde{u}_{22}     &    &  \\
			\vdots                               & 	\vdots                        &          & \ddots & 	\vdots     \\
			(\frac{n}{p})\zeta^{-n}	        & \widetilde{u}_{n1}	                      &     & \ldots&	 \widetilde{u}_{nn}\\
		\end{array} \right),
	\end{displaymath}
	where 
	$$	\widetilde{u}_{ij}=u_{ij}\cdot  \(\frac{i}{p}\)\zeta^{i}\cdot  \(\frac{j}{p}\)\zeta^{j} = \frac{(\frac{j}{p})\zeta^{j}+(\frac{-1}{p})(\frac{i}{p})\zeta^{i}}  { 1+(\frac{-1}{p})(\frac{i}{p})\zeta^{i}(\frac{j}{p})\zeta^{j} }.$$
	Then cancel those $(\frac{-1}{p})$ from the second column  to the last column by the first column we obtain 
	\begin{align}
		\det (GHG) &=\det 	\left( \begin{array}{cccccccc}
			1                                           & 0                & 0    &\ldots   &  0\\
			(\frac{1}{p})\zeta^{-1}	                 &	 \widetilde{u}_{11} -	(\frac{-1}{p})\zeta^{-1}	                      &   \widetilde{u}_{12} -	(\frac{-1}{p})\zeta^{-1}	       &\ldots  &   \widetilde{u}_{1n}-	(\frac{-1}{p})\zeta^{-1}	   \\
			(\frac{2}{p})\zeta^{-2}    &  \widetilde{u}_{21} -	(\frac{-2}{p})\zeta^{-2}                      & \widetilde{u}_{22}  -	(\frac{-2}{p})\zeta^{-2}     &    &    \nonumber \\
			\vdots                               & 	\vdots                        &          & \ddots & 	\vdots     \\
			(\frac{n}{p})\zeta^{-n}	        &\widetilde{u}_{n1}-	(\frac{-n}{p})\zeta^{-n}	 	                      &     & \ldots&	 \widetilde{u}_{nn}-	(\frac{-n}{p})\zeta^{-n}	 \\
		\end{array} \right)    \nonumber  \\
		&= \det \bigg[\widetilde{u}_{ij}-	\(\frac{-i}{p}\)\zeta^{-i}  \bigg]_{1\le i,j\le n} \nonumber\\
		&= \det \bigg[\(\frac{-i}{p}\) \cdot \frac{	 \zeta^{-i}\cdot (\zeta^{2i}-1)}  {1+ (\frac{i}{p})\zeta^{i} \cdot (\frac{-1}{p} )(\frac{j}{p})\zeta^{j} }\bigg]_{1\le i,j\le n}\nonumber\\
		&=\(\frac{n!}{p}\) \zeta^{-n(n+1)/2}\cdot  \prod_{j=1}^{n}(1-\zeta^{2j}) \cdot \det \bigg[ \frac{1}{1+ (\frac{i}{p})\zeta^{i} \cdot (\frac{-1}{p})(\frac{j}{p})\zeta^{j}} \bigg]_{1\le i,j\le n}.      \nonumber \\\nonumber
	\end{align}
	Let 
	$$N=\bigg[ \frac{1}{ 1+ (\frac{i}{p})\zeta^{i} \cdot (\frac{-1}{p})(\frac{j}{p})\zeta^{j}} \bigg]_{1\le i,j\le n}. $$  
	As
	$$\det G= \prod_{i=1}^{n}     \(\frac{i}{p}\) \zeta^{i}, $$
	we have 
	\begin{align}
		\det H&  = \(\frac{n!}{p}\)\zeta^{-3n(n+1)/2} \cdot  \prod_{j=1}^{n}(1-\zeta^{2j}) \cdot \det N. \nonumber
	\end{align}
	
	We next show that  
	$\prod_{j=1}^{n}(1-\zeta^{2j}) \cdot \det N$ is a unit in $\mathbb{Z}[\zeta]$. Applying Lemma \ref{Lem. to prove det N} with $u_i= (\frac{i}{p}) \zeta^{i}$ and $v_j= (\frac{-1}{p})(\frac{j}{p}) \zeta^{j}$ to the matrix $N$, we obtain that $\det N$ is equal to 
	\begin{align}
		&\det \bigg[ \frac{1}{ 1+ (\frac{i}{p})\zeta^{i} \cdot (\frac{-1}{p} )(\frac{j}{p})\zeta^{j}} \bigg]   \nonumber  \\
		=&\prod_{1\le i<j\le n}^{}\(   \(\(\frac{i}{p}\)\zeta^{i}-\(\frac{j}{p}\)\zeta^{j}\)  \(\(\frac{-1}{p}\)\(\frac{j}{p}\) \zeta^{j}-\(\frac{-1}{p}\)\(\frac{i}{p}\)\zeta^{i}\)\)   \nonumber   \\
		& \quad \times\prod_{1\le i,j\le n}^{} \(1+\(\frac{-1}{p}\) \(\frac{i}{p}\)\zeta^{i}\cdot \(\frac{j}{p}\)\zeta^{j}  \)^{-1}   \nonumber \\
		=& (-1)^{(n+1)\binom{n}{2}}\cdot  \prod_{1\le i<j\le n}^{}     \(\(\frac{j}{p}\) \zeta^{j}-\(\frac{i}{p}\)\zeta^{i}\)^{2} \cdot  \prod_{1\le i,j\le n}^{} \(1+\(\frac{-1}{p}\) \(\frac{i}{p}\)\zeta^{i}\cdot \(\frac{j}{p}\)\zeta^{j}  \)^{-1}.  \nonumber
	\end{align}
	Following Vsemirnov's notations, let 
	\begin{equation*}
		f_1= \prod_{1\le i<j\le n}^{}\(   \(\frac{j}{p}\)\zeta^{j}- \(\frac{i}{p}\)\zeta^{i} \)
	\end{equation*}
	and
	\begin{equation*}
		f_2= \prod_{1\le i<j\le n}^{}\( 1+ \(\frac{-1}{p}\) \(\frac{i}{p}\)\zeta^{i}\(\frac{j}{p}\)\zeta^{j} \).
	\end{equation*}
	Then one can verify that 
	\begin{align}
		\prod_{1\le i,j\le n}^{}\( 1+ \(\frac{-1}{p}\) \(\frac{i}{p}\)\zeta^{i}\(\frac{j}{p}\)\zeta^{j} \)=f_2^2\cdot \prod_{j=1}^{n}\(1+\(\frac{-1}{p}\)\zeta^{2j}\), \nonumber
	\end{align}
	and 
	\[  
	\prod_{j=1}^{n}(1-\zeta^{2j}) \cdot \det N=(-1)^{(n+1)\binom{n}{2}}\cdot f_1^2 f_2^{-2}\cdot \prod_{j=1}^{n}\(1+\(\frac{-1}{p}\)\zeta^{2j}\)^{-1}\cdot \prod_{j=1}^{n}\(1-\zeta^{2j}\). 
	\]
	Now we separate our discussion into two cases.	If $p\equiv 3 \pmod 4$, then
	\begin{align}
		&\prod_{j=1}^{n}(1-\zeta^{2j}) \cdot \det N   \nonumber \\
		&=(-1)^{(n+1)\binom{n}{2}}\cdot f_1^2 f_2^{-2}       \nonumber\\
		& =(-1)^{(n+1)\binom{n}{2}}\cdot \bigg(   \prod_{1\le i<j\le n}^{} \frac{ (\frac{j}{p})\zeta^{j}- (\frac{i}{p})\zeta^{i}}{( 1+ (\frac{-1}{p}) (\frac{i}{p})\zeta^{i}(\frac{j}{p})\zeta^{j}}    \bigg)^2  \nonumber \\
		& =(-1)^{(n+1)\binom{n}{2}}\cdot \bigg(
		\prod_{1\le i<j\le n}^{}     \(\frac{j}{p}\) \zeta^{j}  \cdot   \prod_{1\le i<j\le n}^{} \frac{   1- (\frac{i}{p})\zeta^{i}(\frac{j}{p})\zeta^{-j}   }{1- (\frac{i}{p})\zeta^{i}(\frac{j}{p})\zeta^{j} }    \bigg)^2  \nonumber \\
		& =	(-1)^{(n+1)\binom{n}{2}}\cdot		\prod_{1\le i<j\le n}^{}   \zeta^{2j}  \cdot \bigg(   \prod_{1\le i<j\le n}^{} \frac{1- (\frac{i}{p})\zeta^{i}(\frac{j}{p})\zeta^{-j}  }{1- (\frac{i}{p})\zeta^{i}(\frac{j}{p})\zeta^{j} }    \bigg)^2.   \nonumber\\
		\nonumber
	\end{align}
	For any $1\le i<j\le n$, it is clear that 
	
	$$\frac{1- (\frac{i}{p})\zeta^{i}(\frac{j}{p})\zeta^{-j}  }{1- (\frac{i}{p})\zeta^{i}(\frac{j}{p})\zeta^{j} }\in\mathbb{Z}[\zeta]^{\times}.$$
	This implies that $\prod_{j=1}^{n}(1-\zeta^{2j}) \cdot \det N $ is a unit.
	
	In the case $p\equiv 1\pmod 4$, we have $\tau =\sqrt{p}$ and by (\ref{mu})
	\begin{align}
		&\prod_{j=1}^{n}(1-\zeta^{2j}) \cdot \det N   =(-1)^{(n+1)\binom{n}{2}}\cdot \mu\cdot \prod_{j=1}^{n}(1+\zeta^{2j})^{-1} \cdot \sqrt{p} f_1^2 f_2^{-2}.   \nonumber
	\end{align}
 Combining this with  Vsemirnov's result
	$\sqrt{p}f_1^2 f_2^{-2}  \in \mathbb{Z}[\zeta]^\times $, we conclude that $\prod_{j=1}^{n}(1-\zeta^{2j}) \cdot \det N $ is a unit in $\Z[\zeta]$.
	
	In view of the above, we have completed the proof of $\det M_p\in\mathbb{Z}[\zeta]^{\times}$. \qed

	\section{Second part of the proof :  $\det M_p \mod p$}
	In this section we prove the congruence (\ref{mod}).
	We first introduce two lemmas which we shall use in this section.
	\begin{lemma}[The Matrix-Determinant Lemma]
		Let $H$ be an $n\times n$ matrix, and let $\u,\v$ be two $n$-dimensional column vectors. Then
		\begin{equation*}
			\det (H+\u\v^T)=\det H+ \v^T \adj(H)\u,
		\end{equation*}
		where $\adj (H)$ is the adjugate  matrix of $H$ and $\v^T$ is the transpose of $\v$.	
		\label{mdl}
	\end{lemma}
	Readers may refer to \cite{mdl} for a detailed proof of the Matrix-Determinant Lemma. We also need the following result of Sun. 
	\begin{lemma}{\rm (Sun \cite[Theorem 1.3]{Sun-ffa-per} )} Let $p$ be an odd prime. Then the following results hold. 
		
		\noindent	{\rm (i)} 
		If $p\equiv 1 \pmod 4$, then
		\begin{equation}
			\prod_{k=1}^{(p-1)/2}(1-\zeta^{k^2})=\sqrt{p}\varepsilon_p^{-h(p)},
		\end{equation}
		\begin{equation}
			\prod_{1\le j<k\le (p-1)/2}(\zeta^{j^2}-\zeta^{k^2})^2=(-1)^{(p-1)/4}p^{(p-3)/4}\varepsilon_p^{h(p)},
		\end{equation}
		where $\varepsilon_p>1$ and $h(p)$ are the fundamental unit and the class number of the real quadratic field $\mathbb{Q}(\sqrt{p})$ respectively.
		
		\noindent	{\rm (ii)} 	
		If $p\equiv 3\pmod 4$, then
		\begin{equation}
			\prod_{k=1}^{(p-1)/2}(1-\zeta^{k^2})=(-1)^{(h(-p)+1)/2}\sqrt{p}{\bf i},
		\end{equation}
		and
		\begin{equation}
			\prod_{1\le j<k\le (p-1)/2}(\zeta^{j^2}-\zeta^{k^2})^2=(-p)^{(p-3)/4},
		\end{equation}
		where $h(-p)$ is the class number of the imaginary quadratic field $\mathbb{Q}(\sqrt{-p})$.
		\label{sun-lemma-3}
	\end{lemma}
	
	Now we state our proof.
	
	{\noindent\bf Second part of the proof of Theorem \ref{th}.} Recall that $\mathfrak{p}$ is the unique prime ideal of $\mathbb{Z}[\zeta]$ lying over $p \mathbb{Z}$ such that  $\zeta \equiv 1 \pmod {\mathfrak{p}}$ and 
	$\mathfrak{p}\cap \mathbb{Z}=p \mathbb{Z}$.
	
	If $a\equiv 0\pmod p$, then
	\begin{align}
		\frac{1}{\tau}\(1+2\sum_{k=1}^{(p-1)/2}\zeta^{ak^2}\)=\frac{1}{\tau}\(1+2 \cdot \frac{p-1}{2}\)=
		\frac{p}{\tau}=\overline{\tau}.\nonumber
	\end{align}
	Also
	$$
	\overline{\tau}=\sum_{k=1}^{p-1}\left(\frac{k}{p}\right)\zeta^{-k}
	\equiv\sum_{k=1}^{p-1}\left(\frac{k}{p}\right)\cdot 1
	\equiv 0\pmod {\mathfrak{p}}.
	$$
	Thus, by the above and the Galois theory, for any integer $a$ we have 
	\begin{equation}\label{equiv-2}
		\(\frac{a}{p}\)\equiv \frac{1}{\tau}\(1+2\sum_{k=1}^{(p-1)/2}\zeta^{ak^2}\) \pmod {\mathfrak{p}}.
	\end{equation}
	Recall that $n=(p-1)/2$. By (\ref{equiv-2}), for any $1\le i\le n$ and $0\le j\le n$,
	\begin{align}
		\(\frac{i-j}{p}\) &\equiv \frac{1}{\tau}\(1+2\sum_{k=1}^{n}\zeta^{(i-j)\cdot k^2} \) \nonumber\\
		&\equiv \frac{1}{\tau}\(-1+2\sum_{k=0}^{n}\zeta^{(i-j)\cdot k^2} \) \pmod {\mathfrak{p}}.\nonumber
	\end{align}
	Now we define $\widetilde{M_p}=[m_{ij}]_{0\le i,j\le n}$ with
	$m_{1j}=-1$ for $0\le j\le n$ and
	$$m_{ij}=-1+2\sum_{k=0}^{n}\zeta^{(i-j)\cdot k^2} \quad \text{for}\quad i>1.$$
	
	\noindent
	Comparing $\widetilde{M_p}$ with $M_p$, we have 
	\begin{equation}\label{ideal}
		\det M_p\equiv \frac{-1}{\tau^{n}}\cdot \det \widetilde{M_p} \pmod {\mathfrak{p}}.
	\end{equation}
	We next evaluate $\det \widetilde{M_p}$. First note that
	\begin{equation}\label{mptilde}
		\widetilde{M_p}=-{\bm v}\cdot {\bm v}^{T}+2\cdot A\cdot B^T,
	\end{equation}
	where ${\bm v}=(1,1,\ldots,1)^T$ is an $n$-dimensional vector, $B=[\zeta^{-ij^2}]_{0\le i,j\le n}$ and $A$ is the matrix obtained from $[\zeta^{ij^2}]_{0\le i,j\le n}$ via replacing all entries in the first row by zero. Applying Lemma \ref{mdl} to (\ref{mptilde}), we obtain 
	\begin{align*}
		\det \widetilde{M_p}&=\det (2\cdot A\cdot B^T)-{\bm v}^{T} \mathrm{adj}(2\cdot A\cdot B^T) {\bm v} \\
		&=-2^{\frac{p-1}{2}}{\bm v}^{T}\cdot \mathrm{adj}(B)^T \mathrm{adj}(A)\cdot {\bm v}\\
		&=-2^{\frac{p-1}{2}}\cdot \sum_{k=0}^{n}\det B_{(k)} \det A_{(k)},
	\end{align*}
	where $A_{(k)}$ (resp. $B_{(k)}$) is the matrix obtained from $A$ (resp. $B$) via replacing all the entries in the $(k+1)$-th column by $1$. 
	
	Note that all entries in the first column of $B$ are $1$. Hence $\det B_{(k)}=0$ for any $k\ge 1$. It suffices to calculate the determinants of $B_{(0)}=B$ and $A_{(0)}$.
	For convenience, write $\eta_j=\zeta^{j^2}$ for $j=1,2,\ldots,n$. Then
	\begin{equation}
		A_{(0)}=
		\left(
		\begin{array}{ccccc}
			
			1&0&0&\ldots&0\\
			1&\eta_1&\eta_2&\ldots&\eta_n\\
			1&\eta_1^2&\eta_2^2&\ldots&\eta_n^2\\
			\vdots&\vdots&\vdots&\ddots&\vdots\\
			1&\eta_1^n &\eta_2^n&\ldots&\eta_n^n\\
		\end{array}
		\right)\nonumber
	\end{equation}
	and
	\begin{equation}
		B=
		\left(
		\begin{array}{ccccc}
			
			1&1&1&\ldots&1\\
			1&\eta_1^{-1}&\eta_2^{-1}&\ldots&\eta_n^{-1}\\
			1&\eta_1^{-2}&\eta_2^{-2}&\ldots&\eta_n^{-2}\\
			\vdots&\vdots&\vdots&\ddots&\vdots\\
			1&\eta_1^{-n} &\eta_2^{-n}&\ldots&\eta_n^{-n}\\
		\end{array}
		\right).\nonumber
	\end{equation}
	One can verify that 
	\begin{align}
		\det A_{(0)}=\prod_{k=1}^{n}\eta_k\cdot V(\eta_1,\eta_2,\ldots,\eta_n)=\prod_{1\le j<k\le n }(\zeta^{k^2}-\zeta^{j^2}).\nonumber
	\end{align}
	For any matrix $M=[a_{ij}]_{1\le i,j\le r}$, let $\overline{M}=[\overline{a_{ij}}]_{1\le i,j\le r}$, where $\overline{z}$ denotes the complex conjugate of a complex number $z$. Then 
	\begin{align}
		\det B&=V(1,\eta_1^{-1},\ldots,\eta_n^{-1})=\overline{V(1,\eta_1,\eta_2,\ldots,\eta_n)}  \nonumber\\
		&=(-1)^{n}\prod_{k=1}^{n}\overline{(1-\zeta^{k^2})}\cdot \prod_{1\le j<k\le n}\overline{(\zeta^{k^2}-\zeta^{j^2})}.\nonumber
	\end{align}
	This implies 
	\begin{equation}\label{mptidle}
		\det \widetilde{M_p}=-(-2)^{n}\prod_{k=1}^{n}\overline{(1-\zeta^{k^2})}\cdot \bigg|\prod_{1\le j<k\le n}(\zeta^{k^2}-\zeta^{j^2})\bigg|^2.
	\end{equation}
	
	We now divide the remaining proof into two cases. 
	
	{\bf Case I:} $p\equiv 1 \pmod 4$
	
	\noindent
	In this case, it follows from Lemma \ref{sun-lemma-3} that 
	\begin{equation}
		\prod_{k=1}^{n}(1-\zeta^{k^2})=\sqrt{p}\varepsilon_p^{-h(p)}\nonumber
	\end{equation}
	and
	\begin{equation}
		\bigg|\prod_{1\le j<k\le n}(\zeta^{k^2}-\zeta^{j^2})\bigg|^2=p^{(p-3)/4}\varepsilon_p^{h(p)}.\nonumber
	\end{equation}
	Combining (\ref{ideal}) and (\ref{mptidle}) with the above, we obtain 
	\begin{align*}
		\det(M_p)&\equiv \frac{-1}{\tau^{n}}\cdot \det(\widetilde{M_p})\\
		&=\frac{-1}{\tau^{n}}\cdot (-1)\cdot 2^{(p-1)/2}	\cdot p^{(p-1)/4}\\
		&=2^{(p-1)/2}\\
		&\equiv \(\frac{2}{p}\)
		\equiv (-1)^{\frac{p-1}{4}}\pmod {\mathfrak{p}}.
	\end{align*}
	
	\noindent
	As $\mathfrak{p}\cap \mathbb{Z}=p\mathbb{Z}$ and 
	$\det M_p$ and $(-1)^{\frac{p-1}{4}}$ are both integers, we have
	$$\det M_p\equiv (-1)^{\frac{p-1}{4}}\pmod p.$$
	
	\noindent
	{\bf Case II:} $p\equiv 3 \pmod 4$
	
	\noindent
	In this case, by Lemma \ref{sun-lemma-3} it is easy to see that 
	\begin{equation}
		\prod_{k=1}^{n}(1-\zeta^{k^2})=(-1)^{(h(-p)+1)/2}\sqrt{p}{\bf i} \nonumber
	\end{equation}
	and
	\begin{equation}
		\bigg|\prod_{1\le j<k\le n}(\zeta^{k^2}-\zeta^{j^2})\bigg|^2=p^{(p-3)/4}. \nonumber
	\end{equation}
	By the method used in the case $p\equiv 1 \pmod 4$  one can verify that 
	\begin{align*}
		\det M_p &\equiv \frac{-1}{\tau^{n}}\cdot \det \widetilde{M_p}\\
		&=\frac{-1}{\tau^{n}}\cdot 2^{(p-1)/2}	\cdot (-1)^{(h(-p)-1)/2}\sqrt{p}{\bf i}\cdot p^{(p-3)/4}\\
		&\equiv(-1)^{(h(-p)-1)/2}\cdot  \(\frac{2}{p}\)\cdot (-1)^{\frac{p+1}{4}}\\
		&\equiv(-1)^{(h(-p)-1)/2}\pmod {\mathfrak{p}}.
	\end{align*}
	Consequently, $\det M_p \equiv(-1)^{(h(-p)-1)/2}\pmod p$.
	
	In view of the above, we have completed the proof of (\ref{mod}).\qed 
	
	\section{Concluding Remarks}
	
	Sun \cite[Conjecture 4.6(iii)]{S19} also defined $N_p$ as the matrix obtained from $[(\frac{i-j}{p})]_{1\le i,j\le (p-1)/2}$ via replacing all the entries in the first row by $1$. Sun conjectured that 
	$$\det N_p=(-1)^{\lfloor(p+3)/4\rfloor}\det M_p,$$
	where $\lfloor\cdot\rfloor$ is the floor function. 
	
	Although $M_p$ and $N_p$ have very similar structures, it seems that the technique used in this paper cannot be applied to $\det N_p$ directly. We believe that calculating $\det N_p$ is quite challenging.
	\medskip
	
	{\noindent\bf Declaration of competing interest} The authors declare that they have no conflict of interest.
	\medskip
	
	{\noindent\bf Data availability} No data was used for the research described in the article.
	\medskip
	
	\Ack\ The authors would like to thank two referees for helpful comments. We are very grateful to Prof. Zhi-Wei Sun for telling us this interesting conjecture.

\end{document}